\documentclass[a4paper,10pt]{article}
\usepackage{ amstext, amsfonts, amssymb, a4}
\usepackage{amscd, amsthm}
\usepackage{epsfig,epsf}
\usepackage[french]{babel}
\usepackage{epsfig,epsf}
\input{epsf.tex}
\usepackage[latin1]{inputenc}
\usepackage[french]{babel}
\usepackage{amsthm}
\usepackage{enumerate}

\newcommand{\pl}{\mathbb{P}}

\newcommand{\cl}{\mathbb{C}}

\newcommand{\pq}{\nabla}

\newcommand{\rarrow}{\rightarrow}
\newcommand{\lrarrow}{\longrightarrow}
\newcommand{\lac }{\left\{ }
\newcommand{\rac }{\right\} }
\newcommand{ \rp }{\right) }
\newcommand{ \lp }{\left( }
\newcommand{ \cf }{\emph{cf.} }
\newcommand{\ol}{\mathcal{O}}
\newcommand{\dx}{\partial_x}
\newcommand{\dy}{\partial_y}

\newcommand{\W}{\mathcal{W}}

\newcommand{\espace}{\\ $\left.\right.$ \\\indent }

\newcommand{\esp}{\\ \indent }
\newcommand{\vesp}{\vspace{1ex} \esp }
\newcommand{\tir}{\nolinebreak\mbox{-}\nolinebreak}
\begin{document}

\title{Détermination du rang des tissus du plan \\et autres invariants géométriques}
\author{Olivier Ripoll}
\date{\footnotesize Octobre 2004}
\maketitle
\footnotesize\textbf{Résumé :}\\
Soit~$\W(d)$ un~$d\tir$tissu non singulier du plan présenté par une équation différentielle du premier ordre à coeff\-icients dans~$\cl\lac x,y\rac$ de la forme 
\begin{center}
$F(x,y,y'):=a_0(x,y)\cdot (y')^d+a_1(x,y)\cdot (y')^{d-1}+\cdots+a_d(x,y)=0,$
\end{center} avec~$d\geq3$ et dont on note~$(E,\pq)$ sa connexion associée. Nous montrons que la trace de la courbure de~$(E,\pq)$ est la somme des courbures de Blaschke des~$3\tir$tissus extraits. En outre nous indiquons comment la courbure rend compte de la linéarisabilité du tissu~$\W(d)$. Notre résultat principal est un procédé explicite de détermination pour~$d$ quelconque du rang de~$\W(d)$, à partir des coeff\-icients de~$F$. En application, nous retrouvons égalemment des résultats connus en géométrie des tissus.\vesp 
\textbf{Abstract :}\\
Let~$\W(d)$ be a non singular~$d\tir$web in the plane with~$d\geq3$, presented by a first order differential equation of the type $F(x,y,y'):=a_0(x,y)\cdot(y')^d+a_1(x,y)\cdot(y')^{d-1}+\cdots+a_d(x,y)=0,$ where $a_i\in\cl\lac x,y\rac$ and let~$(E,\pq)$ be the connection associated to~$F$. We show that the trace of its curvature is the sum of the Blaschke curvatures of extracted~$3\tir$webs of $\W(d)$. Our main result is an explicit determination of the rank of $\W(d)$. We also recover some well known results in web geometry.\normalsize
\section*{Abridged english version}
Let~$\W(d)$ be a non singular~$d\tir$web in~$(\cl^2,0)$, \emph{presented} by a f\-irst order differential equation with analytic coeff\-icients
$$F(x,y,y'):=a_0(x,y)\cdot (y')^d+a_1(x,y)\cdot (y')^{d-1}+\cdots+a_d(x,y)=0.$$ \indent This \emph{implicit} presentation leads to an explicit study of some invariants of the web~$\W(d)$ as the rank, and f\-its in with the geometric study of differential equations.\esp  Using the connection~$(E,\pq)$ associated to $F$ and introduced in~[4], we can first interpret some of the coeff\-icients of its curvature matrix.\\ \indent We prove for at least~$d=4\,,5\mbox{ and }6$, that the \emph{Blaschke-Chern curvature} of~$\W(d)$, which is the trace of the curvature~$K$ of~$(E,\pq)$, is an invariant of the web, \emph{i.e.} it only depends on the web and not on the presentation chosen. Furthermore, we have the following theorem\vspace{1ex}\\
\bfseries Theorem 1. (Trace formula) \normalfont\itshape The Blaschke-Chern curvature of a~$d\tir$web~$\W(d)$ is the sum of all the Blaschke curvatures of extracted~$3\tir$webs of~$\W(d)$, for at least\\$d=4\,,5\mbox{ and }6$.\normalfont \vesp
Interpreting the coeff\-icients of the connection matrix in any adapted basis, we can show that for a non singular planar~$4\tir$web, there exists an adapted basis where the following expression of the curvature matrix holds~:~$$K=\left(\matrix{ k_1 &\dx(k_1)+L_1 & \dy(k_1)+L_2 \cr0 & 0 & 0 \cr0& 0 & 0\cr}\right)dx\wedge dy$$ where~$L_1=L_2=0$ are the conditions for a~$4\tir$web to be linearisable~(\cf [3]). Hence we can deduce from this the Poincaré theorem~:\emph{ A~$4\tir$web of maximal rank is linearisable}.\vesp In the classical case where two families of leaves of~$\W(4)$ are respectively given by~$x=const.$ and by~$y=const.$, the previous methods still hold and we get an explicit expression of the coeff\-icients of the curvature matrix, involving the Blaschke curvatures of extracted~$3\tir$webs. This allows us to prove once again the trace formula and other results.\vesp In the general case of a $d\tir$web~$\W(d)$, the existence of an adapted basis of~$(E,\pq)$ gives just one integrability relation, leading to the construction of a matrix~$(k_{m\ell})$ of dimension~$\pi_d$ with coeff\-icients in $\cl\lac x,y\rac$ which gives explicitely the rank of~$\W(d)$. Precisely, using the special construction of this matrix, and the fact that~$\mbox{Ker}\,\pq$ is a local sytem, we obtain, thanks to Nakayama's lemma, the following result~:\vspace{1ex}\\
\bfseries Theorem 2. (Determination of the rank of a non singular planar~$d\tir$web) \normalfont\itshape The following equality holds~:~$\overline{K}:=\mbox{Ker}\,(k_{m\ell})=\ol\otimes_\cl\mbox{Ker}\,\pq$, and in particular,~$rank\,\W(d)=corank\,(k_{m\ell})$.\normalfont \vesp
For instance, this result allows us to recover a weak version of a theorem due to G. Bol~: \emph{A hexagonal~$4\tir$web is of maximal rank}.
\section{Introduction à la géométrie des tissus du plan}
\indent La géométrie des tissus est consacrée à l'étude des familles de feuilletages en position générale et à leurs invariants.\\ \indent Soit~$\ol=\cl\{x,y\}$ l'anneau des séries convergentes à coeff\-icients complexes à deux variables. Un~$d\tir$tissu~$\W(d)$ non singulier de~$(\cl^2,0)$ est la donnée d'une famille de \emph{feuilles}, germes de courbes de niveau~$\{F_i(x,y)=cste\}\mbox{ où } F_i\in\ol$ vérif\-ie~$F_i(0)=0$ pour tout~$1\leq i\leq d$, et satisfaisant l'hypothèse de position générale $$dF_i(0)\wedge dF_j(0)\,\neq\,0$$ pour~$1\leq i<j\leq d$. On considère une équation différentielle du premier ordre de la forme~$(1)$ suivante~:
\begin{center}
$F(x,y,y'):=a_0(x,y)\cdot (y')^d+a_1(x,y)\cdot (y')^{d-1}+\cdots+a_d(x,y)=0\qquad (1)$
\end{center}
où~$F(x,y,p)\in\ol[p]$ est sans facteurs multiples et telle que~$a_0\neq0$. \vesp On note~$R\in\ol$ le~$p\tir$résultant de~$F$ avec~$R=(-1)^{\frac{d(d-1)}{2}}\cdot a_0\cdot\Delta$ où~$\Delta$ est son~$p\tir$discriminant. En dehors du lieu~$\lac R=0\rac$, les~$d$ courbes intégrales d'une équation différentielle de la forme~(1) déf\-inissent un~$d\tir$tissu non singulier, en vertu du théorème de Cauchy. Réciproquement, la donnée d'un~$d\tir$tissu non singulier de~$(\cl^2,0)$ permet de construire, quitte à effectuer un changement de variables, une équation différentielle$$F(x,y,y'):=\prod_{i=1}^d(\dy(F_i)y'+\dx(F_i))=0$$ de la forme~(1) vérif\-iant~$R(0)\neq0$ et dont les~$d$~solutions au voisinage de~$0$ ont les pentes $p_i(x,y):=-\dx(F_i)/\dy(F_i)\in\ol$ des feuilles du tissu.\vesp Ainsi, tout tissu du plan est \textit{implicitement présenté} par une équation différentielle de la forme~(1), à un inversible près. \normalfont Cette approche ne privilégie aucune des feuilles du tissu et s'inscrit dans le cadre de l'étude géométrique des équations différentielles. Cette étude sera locale, au voisinage de~$0\in\cl^2$ et l'on suppose \itshape désormais \normalfont que~$R(0)\neq0$.\\\,\\
$\bullet$ \bfseries Invariants des tissus du plan. \normalfont Soit~$\mathcal{A}(d)$ le~$\cl\tir$espace vectoriel des \emph{relations abéliennes} du tissu non singulier~$\W(d)$ liant les normales aux feuilles et à coeff\-icients constants sur celles-ci, déf\-ini par   
$$ \mathcal{A}(d):=\big\{(g_1(F_1),\ldots g_d(F_d))\in\ol^d\mbox{ tel que }(g_i)_{i=1\ldots d}\in\cl\lac t\rac \mbox{ et }\sum_{i=1}^dg_i(F_i)dF_i=0\big\}.$$

\indent On a la majoration classique \emph{optimale} suivante~: 
\begin{center}\vspace{1ex}
$rg\,\W(d):=dim_\cl\mathcal{A}(d) \leq \pi_d:=\frac{1}{2}\,(d-1)(d-2)$
\end{center}
et l'entier~$rg\,\W(d)$ déf\-ini ci-dessus est un \emph{invariant} du tissu~: il ne dépend que de~$\W(d)$ et non du choix des fonctions~$F_i$, ou pareillement, de l'équation différentielle qui le présente. On l'appelle le \emph{rang} du tissu non singulier~$\W(d)$ ; nous verrons qu'en général, le rang d'un tissu est nul.\vesp 
On dit qu'un~$d\tir$tissu~$\mathcal{L}(d)$ est \emph{linéaire} si toutes ses feuilles sont des germes de droites. Un tissu \emph{algébrique}~$\mathcal{L}_C(d)$ est un tissu linéaire donné par dualité par une courbe algébrique réduite~$C$ de~$\pl^2$. Si~$P(s,t)=0$ est une équation aff\-ine de~$C$, on vérif\-ie que dans un système de coordonnées convenables, le tissu~$\mathcal{L}_C(d)$ est présenté par l'équation~$F(x,y,p)=P(y-px,y)=0$. Un théorème d'Abel permet alors de montrer que les tissus algébriques sont de rang maximal~$\pi_d$. \esp Dans~[2], S. S. Chern et P. A. Griff\-iths ont d'ailleurs plus largement illustré le lien existant entre la géométrie des tissus et la géométrie algébrique. Cependant, l'existence des tissus du plan \emph{exceptionnels}, qui sont de rang maximal mais non algébrisables, et dont le premier exemple lié à la relation fonctionnelle du dilogarithme fut donné par G. Bol,  montre que l'on ne peut réduire l'étude des tissus de rang maximal à celle des courbes algébriques planes (\emph{cf.}~[1] et~[4] pour des développements récents). \vspace{1ex} \\
\indent La détermination du rang d'un tissu non singulier quelconque n'était effective que dans le cas d'un~$3\tir$tissu~$\W(3)$. On sait en effet depuis les travaux de W. Blaschke~[1] qu'il existe une~$1\tir$forme~$\gamma$ vérif\-iant la relation $d\omega_i=\gamma\wedge\omega_i$ où~$\omega_i=\rho_idF_i$, avec~$\rho_i\in\ol^*$ et~$\omega_1+\omega_2+\omega_3=0$. La~$2\tir$forme~$d\gamma$ est, contrairement à $\gamma$, un invariant du tissu appelée \emph{courbure de Blaschke} de~$\W(3)$. On a alors l'équivalence suivante~:~$rg\,\W(3)=1$ si et seulement si~$d\gamma=0$.\vspace{1ex} \esp
Nous proposons dans cette note une méthode explicite pour déterminer le rang d'un~$d\tir$tissu du plan non singulier quelconque, ainsi que l'interprétation d'autres invariants du tissu. Notre point de départ sera la connexion associée à un~$d\tir$tissu introduite par A. Hénaut dans~[4] ainsi que les invariants naturels qu'elle engendre.\vesp
$\bullet$ \bfseries Relations abéliennes \emph{via} l'annulation d'une trace, et connexion associée \normalfont(\emph{cf.}~[4]). Soient~$d\geq3$ et~$\W(d)$ un~$d\tir$tissu non singulier présenté par~$F$ de la forme~(1) et~$S=\lac F(x,y,p)=0\rac$ la surface lisse associée, munie de la projection canonique~$\pi:S\lrarrow(\cl^2,0)$ induite par~$(x,y,p)\rarrow(x,y)$. \esp Soit~$(\Omega_S^\bullet,d)$ le complexe de de Rham des formes différentielles sur~$S$. La formule d'interpolation de Lagrange notamment permet de montrer qu'il existe un isomorphisme de~$\cl\tir$espaces vectoriels de~$\mathcal{A}(d)$ dans 
$$\mathfrak{a}_F=\lac\omega=(b_3\cdot p^{d-3}+b_4\cdot p^{d-4}+\cdots+b_d)\cdot\frac{dy-pdx}{\partial_p(F)}\in\pi_\ast(\Omega_S^1),b_i\in\ol\mbox{ et }d\omega=0\rac,$$grâce auquel une relation abélienne est vue comme l'annulation de la trace d'un élément de~$\mathfrak{a}_F$ relativement au revêtement~$\pi$ de degré~$d$.\esp 
Grâce à la fermeture des éléments de~$\mathfrak{a}_F$ et par identif\-ication polynômiale, l'espace~$\mathfrak{a}_F$ est uniquement déterminé par les solutions analytiques~$b_i$ du système différentiel linéaire homogène~$\mathcal{M}(d)$ suivant~:\begin{center}\vspace{1ex} $\displaystyle\mathcal{M}(d)\;\;\left\{\matrix{\dx (b_d) &+&  A_{1,1}\cdot b_3 + \cdots + A_{1,d-2}\cdot b_d&=&0\cr
\dx (b_{d-1})+\dy(b_{d})&+&   A_{2,1}\cdot b_3 + \cdots + A_{2,d-2}\cdot b_d&=&0\cr
  & & & \vdots& \cr 
\dx (b_{3})+\dy(b_{4})&+&  A_{d-2,1}\cdot b_3 + \cdots + A_{d-2,d-2}\cdot b_d&=&0\cr    
\dy(b_3)&+&A_{d-1,1}\cdot b_3 + \cdots + A_{d-1,d-2}\cdot b_d&=&0\cr}\right.$\end{center}\vspace{1ex}
 
\indent Concrètement, un système de Cramer où l'on reconnait le déterminant de Sylvester donnant le~$p\tir$résultant de~$F$ permet de déterminer les coeff\-icients~$A_{ij}$ dans~$\ol\left[1/R\right]$, et un calcul montre qu'ils sont en fait à pôles sur~$\Delta$. Ceci conduit d'ailleurs à une étude des tissus singuliers qui fait l'objet de travaux en cours.\esp  Par la nature de ses symboles, les solutions analytiques de~$\mathcal{M}(d)$ forment un \emph{système local} dont nous sommes amenés à chercher le rang. Pour cela est  construit à partir de~$\mathcal{M}(d)$ grâce à la théorie de Cartan-Spencer un \itshape f\-ibré vectoriel complexe~$E$ de rang~$\pi_d$ sur~$(\cl^2,0)$, \normalfont inclus dans le f\-ibré des jets~$J_{d-2}(\ol^{d-2})$ et une \itshape connexion~$\pq~: E \rarrow\Omega^1\otimes_\ol E$, non nécessairement intégrable, dont les sections horizontales s'identif\-ient à~$\mathfrak{a}_F$. \normalfont De plus, il existe une \emph{base adaptée} de~$(E,\pq)$ telle que sa courbure~$K$ ait pour matrice $$\lp\matrix{k_1 &k_2 & \ldots &k_{\pi_d}\cr 0&0&\ldots&0\cr \vdots &\vdots & &\vdots\cr  0&0&\ldots&0}\rp dx\wedge dy.$$ \indent Le théorème de Cauchy-Kowalevski assure par conséquent que les tissus maximaux sont exactement ceux pour lesquels la connexion associée est intégrable. Dans le cas d'un~$3\tir$tissu, la courbure obtenue est une~$2\tir$forme qui est la \emph{courbure de Blaschke} classique du tissu. Par construction, la connexion~$(E,\pq)$ dépend de la présentation $F$ de $\W(d)$. Cependant, elle permet de déterminer d'autres invariants que le rang, comme nous allons le voir. \vesp Dans toute la suite, on considère un~$d\tir$tissu du plan non singulier~$\W(d)$ présenté par une équation différentielle~$F$ de la forme~(1), et dont on note la connexion associée~$(E,\pq)$ et la courbure~$K$. 
\section{Trace de la matrice de courbure et autres invariants}
Soit~$P_{\W(d)}$ l'unique polynôme à coeff\-icients dans~$\ol$ de degré inférieur ou égal à~$d-1$ vérif\-iant pour tout~$1\leq i\leq d$ l'équation~$P_{\W(d)}(x,y,p_i)=\dx(p_i)+p_i\dy(p_i)$. Les feuilles du tissu sont solutions de l'équation différentielle~$y''=P_{\W(d)}(x,y,y')$ et~$\W(d)$ est linéaire si et seulement si~$P_{\W(d)}$ est nul. Dans~[3] est établi que~$\W(d)$ est linéarisable si et seulement si~$\deg(P_{\W(d)})\leq 3$ et ses coeff\-icients vérif\-ient un système différentiel~$L_1=0$ et~$L_2=0$. Essentiellement par le calcul, on montre la proposition suivante~: \vspace{1ex} \\
\textbf{Proposition 1. }\\
	(\textit{a}) \textit{Le polynôme~$P_{\W(d)}$ est un invariant du tissu~$\W(d)$} ;\\
	(\textit{b}) \textit{Les éléments~$A_{ij}$ du système~$\mathcal{M}(d)$ s'expriment explicitement à partir des coeff\-icients de~$P_{\W(d)}$, et de ceux d'une forme de Pfaff dont la différentielle~$\kappa \,dx\wedge dy$ est un invariant du tissu~$\W(d)$.}
\normalfont\vesp
Par conséquent, on précise les résultats obtenus dans~[4] concernant les \emph{invariants} engendrés par le système~$\mathcal{M}(d)$, et un calcul montre pour au moins~$d=4$ et~$5$, que les coeff\-icients de la matrice de courbure sont des invariants du tissu.\vesp
L'écriture des coeff\-icients~$A_{ij}$ en fonction notamment des coeff\-icients du polynôme $P_{\W(d)}$, permet aussi de retrouver, à la suite de I. Nakai (\emph{cf.}~[5]), pour un~$4\tir$tissu $\W(4)$ du plan non singulier, l'\textit{équivalence} des trois propriétés suivantes~:
\begin{enumerate}
	\item \textit{Pour~$1\leq i \leq 4$, les~$1\tir$formes ~$dF_i$ définissant~$\W(4)$ à un inversible de~$\ol$ près, vérif\-ient la relation $dF_i=s_idF_1+t_idF_4$ avec~$s_i$ et~$t_i$ distincts dans~$\pl^1$} ;
	\item \textit{Les courbures de Blaschke des~$3\tir$tissus extraits de $\W(4)$ sont égales} ;
	\item \textit{ Le birapport des quatre pentes~$p_i$ de $\W(4)$ est constant}.
\end{enumerate}

\indent Par ailleurs, on montre que de tels tissus ne sont \textit{pas} linéarisables, sauf s'ils sont hexagonaux (\emph{i.e.} les quatre~$3\tir$tissus extraits sont de rang~$1$)\vesp
$\bullet$ \bfseries Linéarisation des~$4\tir$tissus. \normalfont L'introduction du polynôme~$P_{\W(d)}$ permet aussi dans le cas où~$d=4$ au moins, d'identif\-ier plus précisément la courbure~$K$, en  montrant la proposition de linéarisation suivante~:\vesp
\textbf{Proposition 2. }\itshape Avec les notations précédentes et pour un~$4\tir$tissu $\W(4)$, il existe une base adaptée de~$(E,\pq)$ telle que la matrice de courbure dans cette base admette la forme suivante \upshape:\itshape $$\left(\matrix{ 3\kappa &3\dx(\kappa)+L_1 & 3\dy(\kappa)+L_2 \cr0 & 0 & 0 \cr0& 0 & 0\cr}\right)dx\wedge dy.$$ 
\normalfont
Par conséquent, et grâce aux propriétés générales de~$P_{\W(d)}$ rappelées ci-dessus,\begin{center} $\W(4)$ \textit{est linéarisable si et seulement si, dans cette base},\\$k_2=3\dx(\kappa)$ et~$k_3=3\dy(\kappa)$.
\end{center}

\indent En fait, on montre dans le cas général que la matrice de courbure d'un tissu linéaire ne dépend que de~$\kappa$ et de ses dérivées partielles. 
Ainsi, on redémontre avec des méthodes propres aux systèmes différentiels, le théorème de Poincaré suivant~: \textit{Un~$4\tir$tissu du plan non singulier de rang maximal est linéarisable.}\vspace{1ex}\\$\bullet$ \bfseries Tissus rectif\-iés. \normalfont Traditionnellement, on peut préférer à la présentation implicite d'un tissu, la donnée d'un \emph{tissu rectif\-ié}, déf\-ini par les courbes de niveaux $$(x,y,F_3(x,y),\ldots,F_d(x,y)).$$

\indent Dans le cas non singulier, un changement de variables permet de se ramener à ce cas et les méthodes précédentes s'y appliquent et donnent le résultat suivant~: \vspace{1ex}\\
\textbf{Proposition 3. }\itshape Pour un~$4\tir$tissu rectif\-ié~$\W(4)$, il existe une base adaptée de~$(E,\pq)$ telle que les coeff\-icients de la matrice de courbure soient  \upshape: $$k_1=K_x+K_y+K_3+K_4\;;\;k_2=\dx(K_y)+v_3K_y\;;\;k_3=\dy(K_x)-v_2K_x$$
\itshape avec~$P_{\W(4)}=-v_1p^3-v_2p^2-v_3p-v_4$ et où~$K_x,K_y,K_3$ et~$K_4$ désignent les courbures de Blaschke des~$3\tir$tissus extraits de~$\W(4)$ indexées par la feuille retirée.\normalfont\vesp  En général, pour $d\geq5$, cette simplicité n'est plus et les deux approches présentées se valent alors en complexité. On remarque dans la proposition précédente que la trace de~$K$ est la somme des courbures de Blaschke des~$3\tir$tissus extraits, ce que l'on va généraliser.\vesp
$\bullet$ \bfseries Sur la trace de la courbure de~$(E,\pq)$. \normalfont En germe dans [4], et de façon analogue à la théorie de Chern-Weil, la trace~$tr(K)$ de la courbure~$K$ se révèle d'une importance propre pour la géométrie du tissu~$\W(d)$ et non pas seulement du f\-ibré~$E$. Comme nous l'avons montré pour le moment par un calcul pour au moins~$d=4,5\mbox{ et }6$, cette trace~$tr(K)$ est un invariant de~$\W(d)$ qu'on appellera la \emph{courbure de Blaschke-Chern} de~$\W(d)$, qui de plus, rend compte des~$3\tir$tissus extraits. En effet, on a le résultat suivant~:\vspace{1ex}\\
\textbf{Théorème 1. (Formule de la trace)}\itshape$\left.\right.$ \\\indent 
Soit~$\W(d)$ un d-tissu non singulier du plan implicitement présenté par~$F$, dont on note~$(E,\pq)$ la connexion associée, et~$K$ sa courbure. On a alors l'égalité suivante  \upshape:\itshape 
\begin{center}
~$\displaystyle tr(K)=\sum_{k=1}^{(^d_3)} d\gamma_k$,
\end{center}
\itshape pour au moins~$d=4,5\mbox{ et }6$, où les~$d\gamma_k$ désignent les courbures de Blaschke des~$3\tir$tissus extraits de~$\W(d)$.\normalfont\vesp
De façon équivalente, le f\-ibré déterminant~$(\det(E),\det(\pq))$ associé à $(E,\pq)$ est donc isomorphe au f\-ibré tensoriel des f\-ibrés en droites~$(\otimes L_i, \otimes \pq_i)$  associés aux~$3\tir$tissus extraits. 
\section{Détermination du rang des tissus du plan~: cas général}
\indent L'existence de bases adaptées de~$(E,\pq)$ permet de donner un moyen effectif de détermination du rang. \esp On note~$\gamma=\gamma_xdx+\gamma_ydy$ la matrice de connexion du tissu dans une base adaptée donnée et on considère la première ligne de la matrice de courbure dans cette base notée $k=(k_1, k_2, \ldots, k_{\pi_d})$. Les sections horizontales de la connexion~$\pq$ notées$$f=\,^t(f_1, f_2, \ldots ,f_{\pi_d})\in E^\pq:=\mbox{Ker}\,\pq$$qui s'identif\-ient aux relations abéliennes du tissu, \emph{via} l'espace~$\mathfrak{a}_F$, vérif\-ient le système différentiel~$df+\gamma f=0.$ La condition d'intégrabilité de ce système est donnée, grâce à la base adaptée, par la \emph{seule} relation $$k\cdot f=k_1f_1+k_2f_2+\cdots+k_{\pi_d}f_{\pi_d}=0\qquad (2).$$ \indent Ces considérations vont nous permettre de démontrer le théorème suivant~: \vspace{1ex}\\   
\textbf{Théorème 2. (Détermination du rang d'un tissu du plan)}\itshape~$\left.\right.$ \\\indent
Soit~$\W(d)$ un d-tissu du plan non singulier présenté par~$F$, dont on note~$(E,\pq)$ la connexion associée, munie d'une base adaptée. Il existe un f\-ibré vectoriel~$\overline{K}$ de rang~$rg\,\W(d)$, noyau d'un endomorphisme de~$\ol^{\pi_d}$, \emph{explicite}, tel que l'on ait $$\overline{K}=\ol\otimes_\cl E^\pq.$$
\indent En particulier, si l'on note~$(k_{m\ell})$ la matrice de cet endomorphisme, le rang du tissu est explicitement donné par l'égalité suivante \upshape: \itshape$rg\,\W(d)=corang\,(k_{m\ell}).$
\normalfont\vspace{1ex}\\
\textit{Démonstration. } Elle utilise le lemme de Nakayama et fondamentalement que~$E^\pq$ est un système local. On considère les~$\pi_d$ équations obtenues de la dérivation jusqu'à l'ordre~$d-3$ de l'équation (2) où l'on pose~$df=-\gamma f$. On obtient ainsi une matrice carrée~$(k_{m\ell})$ d'ordre~$\pi_d$ à coeff\-icients dans $\ol$, dont la première ligne est~$k$. Le~$\ol\tir$module~$\overline{K}:=\mbox{Ker}\,(k_{m\ell})$ est de type f\-ini et par construction, l'inclusion~$\ol\otimes_\cl E^\pq\subseteq \overline{K}$ est vérif\-iée. Pour montrer la relation annoncée, il suff\-it de montrer que l'on a l'égalité  $$\overline{K}=\ol\otimes_\cl E^\pq+\mathfrak{m}\cdot\overline{K}$$ en notant~$\mathfrak{m}$ l'idéal maximal de~$\ol$. On considère des générateurs~$(g_1,\ldots,g_r)$ de~$\overline{K}$ tels que les~$g_i(0)$ soient $\cl\tir$linéairement indépendants. Soit~$g\in\overline{K}$ ; le système précédent de générateurs permet de montrer que si~$g(0)=0$, alors~$g\in\mathfrak{m}\overline{K}$.\\ \indent A l'aide du théorème de Cauchy et grâce à la construction particulière de la matrice $(k_{m\ell})$ issue de la dérivation de la condition d'intégrabilité $(2)$, on montre que si~$g(0)\neq0$, il existe un germe de fonction analytique~$f$ dans~$E^\pq$ telle que\\$f(0)=g(0)$. Ainsi,~$g-f\in\overline{K}$ et vérif\-ie~$(g-f)(0)=0$. Comme précédemment, le système de générateurs permet alors de montrer que~$g-f\in\mathfrak{m}\overline{K}$, ce qui prouve l'égalité voulue entre~$\ol\tir$modules et donne une méthode explicite d'évaluation du rang.~$\Box$ \vesp
Par conséquent, un~$d\tir$tissu non singulier~$\W(d)$ du plan est de rang au moins~$1$ si et seulement si le déterminant de~$(k_{m\ell})$ est nul. Ce \itshape déterminant \normalfont est certainement un invariant du tissu~$\W(d)$, comme cela a été montré dans le cas où~$d=4$.\esp Pour un~$4\tir$tissu et avec les notations précédentes, on vérif\-ie que la matrice~$(k_{m\ell})$ est de la forme $$(k_{m\ell})\;=\; \lp\matrix{k \cr \dx(k)-k\cdot\gamma_x\cr \dy(k)-k\cdot\gamma_y}\rp$$
\indent En écrivant explicitement la matrice~$(k_{m\ell})$ dans ce cas, les résultats précédents permettent de montrer que \emph{le rang d'un $4\tir$tissu~$\W(4)$ dont les courbures de Blaschke des~$3\tir$tissus extraits sont égales et non nulles vérif\-ie l'inégalité~$rg\,\W(4)\leq1$}. On montre également qu'un~$4\tir$tissu de courbure de Blaschke-Chern nulle est de rang distinct de~$2$. Ceci nous permet par exemple de retrouver le résultat suivant:  \vspace{1ex}\\
\textbf{Proposition 4. }(G. Bol, cas~$d=4$)  \textit{Un~$4\tir$tissu hexagonal~$\mathcal{H}(4)$ du plan est de rang maximal}.
\vspace{1ex}\\
\textit{Démonstration. }Un tissu hexagonal~$\mathcal{H}(4)$ admet quatre relations abéliennes à trois termes, dont deux sont nécessairement linéairement indépendantes donc $rg\;\mathcal{H}(4)$ vaut au moins $2$. Or, sa courbure de Blaschke-Chern est nulle en vertu de la formule de la trace, donc~$\mathcal{H}(4)$ ne peut être de rang~$2$, d'après la remarque précédente et par conséquent~$\mathcal{H}(4)$ est de rang maximal.~$\Box$\espace
\noindent\textbf{Références}\small\vspace{1ex}\\
$\left[1\right]$\quad W. BLASCHKE und G. BOL, Geometrie der Gewebe, Springer-Verlag, Berlin, 1938.\\
$\left[2\right]$\quad S. S. CHERN and P. A. GRIFFITHS, Abel's Theorem and Webs, Jahresber. Deutsch. Math.-Verein. \textbf{80} (1978), 13-110.\\
$\left[3\right]$\quad A. H\'ENAUT, Sur la linéarisation des tissus de~$\cl^2$, Topology \textbf{32} (1993), 531-542.\\
$\left[4\right]$\quad A. H\'ENAUT, On planar web geometry through abelian relations and connections, Ann. of Math. \textbf{159} (2004), 425-445.\\
$\left[5\right]$\quad I. NAKAI, Curvature of curvilinear 4-webs and pencils of one forms~: Variation on a theorem of Poincaré, Mayrhofer and Reidemeister, Comment. Math. Helv. \textbf{73} (1998), 177-205.\normalsize
\\\,\\
\footnotesize \noindent Olivier Ripoll\\ LaBAG, UMR 5467\\ 351, cours de la Libération, 33405 Talence, France\\
\footnotesize\textit{Adresse e-mail}~: Olivier.Ripoll@math.u-bordeaux1.fr\normalsize
 \end{document}